\documentclass[a4paper]{article}
\usepackage{amssymb,amsmath,amsthm,verbatim,amscd,enumerate,url}

\DeclareMathOperator{\Supp}{Supp}

\DeclareMathOperator{\Aut}{Aut}
\DeclareMathOperator{\diam}{diam}
\DeclareMathOperator{\girth}{girth}

\title{Box spaces, group extensions and coarse embeddings into Hilbert space}
\author{A. Khukhro\\ \texttt{ak11g08@soton.ac.uk}}

\parskip\bigskipamount

\theoremstyle{plain}
\newtheorem*{prop}{Proposition}

\newtheorem*{thm}{Theorem}

\newtheorem{Thm}{Theorem}
\newtheorem{Lemma}[Thm]{Lemma}
\newtheorem{Prop}[Thm]{Proposition}
\newtheorem{Corollary}[Thm]{Corollary}

\theoremstyle{definition}

\newtheorem*{definition}{Definition}

\newtheorem{Remark}[Thm]{Remark}
\newtheorem{Example}[Thm]{Example}

\begin{document}
\maketitle

\begin{abstract}
We investigate how coarse embeddability of box spaces into Hilbert space behaves under group extensions. In particular, we prove a result which implies that a semidirect product of a finitely generated free group by a finitely generated residually finite amenable group has a box space which coarsely embeds into Hilbert space. This provides a new class of examples of metric spaces with bounded geometry which coarsely embed into Hilbert space but do not have property A, generalising the example of Arzhantseva, Guentner and Spakula.
\end{abstract}


\section*{Introduction}

In this paper, we will mainly be concerned with generalising the example of a bounded geometry metric space which coarsely embeds into Hilbert space but does not have property A, given by Arzhantseva, Guentner and Spakula in \cite{AGS}. In order to do this, we give some basic properties of box spaces, and prove a technical result about extensions. We use this to prove the following theorem.

\begin{thm}
Let $\Gamma= H\rtimes G$ be a finitely generated semidirect product of two residually finite groups $H$ and $G$. If $G$ is amenable and $H$ has a nested sequence of finite index characteristic subgroups with trivial intersection such that the corresponding box space embeds coarsely into Hilbert space, then $\Gamma$ has a box space which coarsely embeds into Hilbert space.
\end{thm}

In particular, the above result applies to \emph{(finitely generated free)-by-cyclic} groups, a large class which includes the fundamental groups of certain 3-manifolds. 

\begin{definition}
Let $(X,d_X)$ and $(Y,d_Y)$ be metric spaces.
$X$ is \emph{coarsely embeddable} in $Y$ if there is a map $F:X\longrightarrow Y$ such that there exist non-decreasing functions $\rho_\pm: \mathbb{R}_+ \longrightarrow \mathbb{R}_+$ with $\lim_{t\rightarrow \infty}\rho_\pm(t)=\infty$ and
\begin{equation*}
\rho_-(d_X(x,x'))\leq d_Y(F(x),F(x'))\leq \rho_+(d_X(x,x'))
\end{equation*}
for all $x,x' \in X$. 
\end{definition}
In the context of coarse geometry, coarse embedding (sometimes also referred to as uniform embedding) is the natural notion of inclusion of one space into another. 

We are particularly interested in discrete metric spaces with bounded geometry which admit a coarse embedding into Hilbert space because by a result of Yu \cite{Yu} the coarse Baum--Connes conjecture holds for such spaces. In \cite{Yu}, Yu introduces \emph{property A}, a non-equivariant notion of amenability, which guarantees coarse embeddability in Hilbert space.

Here we will not give Yu's original definition, but an equivalent definition for spaces with bounded geometry from \cite{Tu}.
\begin{definition}
A discrete metric space $(X,d)$ with bounded geometry has \emph{property A} if and only if for every $R>0$ and $\varepsilon >0$ there exists an $S>0$ and a  function $\phi: X \longrightarrow \ell^2(X)$ such that $\|\phi(x)\|=1$ for all $x \in X$ and such that for all $x_1, x_2 \in X$:
\begin{description}
\item[$(1)$]
if $d(x_1,x_2)\leq R$ then $|1-\left\langle \phi(x_1), \phi(x_2)\right\rangle| \leq \varepsilon$, and
\item[$(2)$]
$\Supp\phi(x) \subset B_{S}(x)$ for all $x\in X$.
\end{description}
\end{definition}

Property A is designed to provide us with a coarse embedding in Hilbert space, and so it is natural to ask whether the converse holds. A counterexample was first given by Nowak \cite{Now}, although his disjoint union of $n$-dimensional cubes $\{0,1\}^n$ over all $n\in \mathbb{N}$ is not of bounded geometry. The question of whether there exist bounded geometry metric spaces without property A which coarsely embed into Hilbert space has recently been answered affirmatively by Arzhantseva, Guentner and Spakula \cite{AGS}. Their elegant example is a carefully chosen \emph{box space} of the free group on two generators.

Let $G$ be a finitely generated residually finite group.
Let $\{K_i\}$ be some collection of finite index subgroups of $G$, for which the intersection $\cap_{n\in \mathbb{N}} K_n$ is trivial.
We will be particularly interested in the case when the $K_i$ are nested normal subgroups of $G$:
\begin{equation*}
G= K_1 \triangleright K_2 \triangleright K_3 \triangleright K_4 \triangleright...
\end{equation*}
Note that given a finitely generated residually finite group, we can always choose a collection $\{K_i\}$ of finite index subgroups with trivial intersection such that the $K_i$ are nested normal subgroups of $G$.
\begin{definition}
The \emph{box space} of $G$ corresponding to $\{K_i\}$, denoted by $\Box_{\{K_i\}} G$, is the disjoint union $\sqcup_i G/K_i$ of finite quotient groups of $G$, where each quotient is endowed with the metric induced by the image of the generating set of $G$, and the distances between the identity elements of two successive quotients are chosen to be greater than the maximum of their diameters. 
\end{definition}
We should think of the box space as a sequence of finite Cayley graphs which are ``strung" onto a thread through the identity elements, so that the distance between two elements from two different Cayley graphs, $\alpha \in G/K_i$ and $\beta \in G/K_j$, is given by $$d(\alpha, \beta)= d_i(\alpha,e_i) + d(e_i,e_j) + d_j(e_j,\beta)=d_i(\alpha,e_i)+\sum_{k=i}^{j-1} d(e_k,e_{k+1})+d_j(e_j,\beta),$$ where $e_i$ denotes the identity element of $K_i$, and $d_i$ denotes the metric on the quotient $G/K_i$. In particular, the distance $d(e_i,e_j)$ is given by $\sum_{k=i}^{j-1} d(e_k,e_{k+1})$.

Note that however we choose the distances between the quotients, we will obtain coarsely equivalent spaces as long as we choose the distances between two successive quotients to be greater than the maximum of their diameters. The properties of $\Box_{\{K_i\}} G$ can vary greatly depending on the choice of $\{K_i\}$. For example, Arzhantseva, Guentner and Spakula's chosen sequence of subgroups gives a box space of $\mathbb{F}_2$ which coarsely embeds in Hilbert space. However, the full box space of $\mathbb{F}_2$ corresponding to the collection of \emph{all} finite index normal subgroups does not, since the box space of $SL(3,\mathbb{Z})$, being a quotient of $\mathbb{F}_2$, coarsely embeds in the full box space $\Box_{\text{all}}\mathbb{F}_2$ of $\mathbb{F}_2$. Thus $\Box_{\text{all}} \mathbb{F}_2$ contains a coarsely embedded expander, for reasons explained below.

There are many remarkable results which link geometric properties of $\Box G$ to analytic properties of $G$. In order to list some of them, we will need the following definitions.

\begin{definition}
A second countable locally compact group $G$ has the \emph{Haagerup property} if there exists a continuous isometric metrically proper action of $G$ on some affine Hilbert space.
\end{definition}
The above is not the original definition of the Haagerup property, but an equivalent property called \emph{a-T-menability}, a term due to Gromov, so called because the Haagerup property is a notion of weak amenability and a strong negation of \emph{Kazhdan's property (T)}. For more on the Haagerup property, see \cite{CJV}. We give the following formulation of property (T), rather than the original definition, since it illuminates the connection between the two properties.

\begin{definition}
A second countable locally compact group $G$ has \emph{property (T)} if every continuous isometric action of $G$ on an affine Hilbert space has a fixed point.
\end{definition}

An exposition of the following can be found in \cite{Roe}.

\begin{prop}
Let $G$ be a finitely generated residually finite group, and let $\{K_i\}$ be a collection of finite index subgroups with trivial intersection. Then:
\begin{itemize}
\item
If $\Box_{\{K_i\}} G$ is coarsely embeddable into Hilbert space, then $G$ has the Haagerup property;
\item
(Guentner) $\Box_{\{K_i\}} G$ has property A if and only if $G$ is amenable.
\end{itemize}
\end{prop}

The first explicit construction of expander graphs by Margulis \cite{Mar} can be rephrased in the following way, using box spaces.
\begin{prop}
Let $G$ be a finitely generated residually finite group with property (T), and equip it with the word length arising from a symmetric generating set $S$ which does not contain the identity. If $\{K_i\}$ is a sequence of nested finite index normal subgroups with trivial intersection, then $\Box_{\{K_i\}} G$ is an expander.
\end{prop}
For instance, a box space of $SL(3,\mathbb{Z})$ is an expander. 

We can see that there is a spectrum of analytic properties of groups which roughly corresponds to geometric properties of box spaces. To summarise, for a residually finite group $G$, and $\{K_i\}$ a sequence of nested finite index normal subgroups of $G$ with trivial intersection, we have
\begin{align*}
G \text{ amenable } &\Longleftrightarrow \Box_{\{K_i\}}G \text{ property A,}\\
G \text{ Haagerup } &\Longleftarrow \Box_{\{K_i\}}G \text{ coarsely embeddable into Hilbert space,}\\
G \text{ property (T) } &\Longrightarrow \Box_{\{K_i\}}G \text{ expander.}
\end{align*}
Note that the last two implications are not reversible.

Another relevant property is \emph{property $(\tau)$}, a weak version of property (T). However, we will not discuss this here, and instead refer the reader to \cite{LZ} for more details.

In what follows, we will say that a metric space is \emph{embeddable} if it embeds coarsely into Hilbert space.

\section*{The construction of Arzhantseva, Guentner and Spakula}

The results in this paper are inspired by the construction of a box space of the free group on two generators which coarsely embeds into Hilbert space \cite{AGS}. We will make use of this result repeatedly.

In \cite{AGS}, the authors define a sequence $\{N_i\}$ of normal subgroups of the free group $\mathbb{F}_2$ inductively: let $N_0:= \mathbb{F}_2$ and let $N_{i+1}:=N_i^2$, where $N_i^2$ denotes the subgroup of $N_i$ generated by all the squares of the elements of $N_i$. Note that each $N_{i+1}$ is a \emph{verbal} subgroup of the previous subgroup $N_i$, and is thus fully invariant and characteristic (in both $N_i$, and the whole of $\mathbb{F}_2$, since being fully invariant and characteristic are transitive properties). In fact, we note that each $N_i$ is a verbal subgroup of $\mathbb{F}_2$, since fully invariant and verbal are equivalent for free groups (see for example 2.3.1 of \cite{Rob}).

By a theorem of Levi (see Proposition 3.3 in Chapter 1 of \cite{LS}), the intersection $\cap N_i$ of all the $N_i$ is trivial. The Cayley graph of $\mathbb{F}_2/N_{i+1}$ is the \emph{$\mathbb{Z}/2\mathbb{Z}$-homology cover} of $\mathbb{F}_2/N_i$, so when we take the quotient of the fundamental group of the Cayley graph of $\mathbb{F}_2/N_i$ (which is isomorphic to $N_i$) by the fundamental group of the Cayley graph of $\mathbb{F}_2/N_{i+1}$ (which is isomorphic to $N_{i+1}$), we obtain a finite direct sum of $\mathbb{Z}/2\mathbb{Z}$'s. 

In \cite{AGS}, Arzhantseva, Guentner and Spakula consider the general setting of a pair of graphs $(\widetilde{X},X)$ where $X$ is a finite graph such that removing a single edge does not disconnect the graph (2\emph{-connected}), and $\widetilde{X}$ is the $\mathbb{Z}/2\mathbb{Z}$-homology cover of $X$. They use the cover $\widetilde{X}$ to induce a wall structure on $X$, and show that for distances smaller than the \emph{girth} of the graph (i.e. the length of the smallest loop), the \emph{wall metric} induced by the wall structure agrees with the natural graph metric. 

In the case of $\mathbb{F}_2$, the box space $\Box_{\{N_i\}}\mathbb{F}_2$ with metric $d_W$ (defined to be the wall metric on each component and the distances between successive quotients taken to be larger than both of their diameters), is coarsely embeddable into Hilbert space since $d_W$ is an effective symmetric normalized negative type kernel on $\Box_{\{N_i\}}\mathbb{F}_2$, in the sense of \cite{Roe}, Chapter 11. Then, since the girth of the graphs $\mathbb{F}_2/N_i$ tends to infinity, the authors conclude that the metric $d$ on $\Box_{\{N_i\}}\mathbb{F}_2$ induced by the word metric of $\mathbb{F}_2$ with its natural generating set is coarsely equivalent to the metric $d_W$, and hence $(\Box_{\{N_i\}}\mathbb{F}_2, d)$ embeds coarsely into Hilbert space. We now make a remark about the generality of this construction, which follows directly from the results of \cite{AGS}.

\begin{Remark}\label{GenAGS}
Given a sequence $\{\widetilde{X}_i\}$ of graphs with increasing diameters such that each $\widetilde{X}_i$ is the $\mathbb{Z}/2\mathbb{Z}$-homology cover of a finite $2$-connected graph $X_i$, one can induce a wall structure and hence a wall metric on each $\widetilde{X}_i$ from $X_i$. The disjoint union $\sqcup \widetilde{X}_i$, metrized using the wall metrics, is coarsely embeddable into Hilbert space. So provided that $\girth(X_i)$ tends to infinity, $\sqcup \widetilde{X}_i$ metrized in the usual way using the natural graph metrics will be coarsely equivalent to $\sqcup \widetilde{X}_i$ with the wall metrics and will thus coarsely embed into Hilbert space. 
\end{Remark}

\section*{Basics}

Recall that given a box space, we obtain a coarsely equivalent space if we change the distances between the quotients, as long as the distances between two successive quotients are chosen to be greater than the maximum of their diameters. Thus, for a finitely generated residually finite group $H$ with a proper left-invariant metric $d$ and a nested sequence of finite index normal subgroups $\{K_i\}$ with trivial intersection, we will write $(\Box_{\{K_i\}}H, d')$ for the metric space obtained by taking the metric induced by $d$ on each finite quotient, with some valid choice of distances between quotients.

\begin{Prop}\label{Metric}
Suppose $d_1$ and $d_2$ are two proper left-invariant metrics on a finitely generated residually finite group $H$. Consider a box space $\Box_{\{K_i\}}H$ of $H$ with respect to some nested sequence of finite index normal subgroups $\{K_i\}$ with trivial intersection. The metric space $(\Box_{\{K_i\}}H, d'_1)$ is coarsely equivalent to $(\Box_{\{K_i\}}H,d'_2)$.
\end{Prop}
\begin{proof}
By Proposition 2.3 of \cite{DG}, the spaces $(H,d_1)$ and $(H,d_2)$ are coarsely equivalent. In other words, the identity mapping from  $(H,d_1)$ to $(H,d_2)$ is a coarse embedding, so there exist non-decreasing functions $\rho_\pm: \mathbb{R}_+ \longrightarrow \mathbb{R}_+$ with $\lim_{t\rightarrow \infty}\rho_\pm(t)=\infty$ such that
\begin{equation*}
\rho_-(d_1(g,h))\leq d_2(g,h)\leq \rho_+(d_1(g,h))
\end{equation*}
for all $g,h \in H$. We need to show that the identity mapping from $(\Box_{\{K_i\}}H, d'_1)$ to $(\Box_{\{K_i\}}H, d'_2)$ is a coarse embedding. 

We can choose the distance between the identity elements of successive quotients to be the same in both spaces, ensuring that this distance is greater than the maximum of the diameters of the two quotients, with respect to both metrics $d'_1$ and $d'_2$.

If two elements $gK_i, hK_j \in \Box_{\{K_i\}}H$ lie in different quotients, we have 
\begin{align*}
d'_1(gK_i,hK_j)&=d'_1(gK_i,eK_i)+d'_1(eK_i,eK_j)+d'_1(eK_j,hK_j)\\
&=d'_1(gK_i,eK_i)+d'_2(eK_i,eK_j)+d'_1(eK_j,hK_j)\\
&\leq 3d'_2(eK_i,eK_j) \leq3d'_2(gK_i,hK_j),
\end{align*}
and vice versa. Thus we have
$$\frac{1}{3}d'_1(gK_i,hK_j)\leq d'_2(gK_i,hK_j) \leq 3 d'_1(gK_i,hK_j).$$

When $gK_i$ and $hK_i$ lie in the same quotient $H/K_i$, choose $k_2$ in $H$ such that its image $k_2K_i$ in $H/K_i$ is equal to $g^{-1}hK_i$ and such that 
\begin{align*}
d'_2(gK_i,hK_i)&=d'_2(eK_i,g^{-1}hK_i):=\min_{b\in H} \{d_2(e,b): bK_i=g^{-1}hK_i \in H/K_i \}\\
&=d_2(e,k_2),
\end{align*}
i.e. the above minimum is attained at $k_2\in H$. Similarly, choose $k_1$ in $H$ such that $k_1K_i=g^{-1}hK_i$ and $$\min_{b\in H} \{d_1(e,b): bK_i=g^{-1}hK_i \in H/K_i \}=d_1(e,k_1).$$
We then have
\begin{align*}
\rho_-(d'_1(gK_i,hK_i))&=\rho_-(d'_1(eK_i,g^{-1}hK_i))=\rho_-(d'_1(eK_i,k_2K_i))\\
&\leq \rho_-(d_1(e,k_2))\leq d_2(e,k_2)\\
&=d'_2(eK_i,g^{-1}hK_i)=d'_2(gK_i,hK_i)
\end{align*}
and
\begin{align*}
d'_2(gK_i,hK_i)&=d'_2(eK_i,g^{-1}hK_i)=d'_2(eK_i,k_1K_i)\\
&\leq d_2(e,k_1)\leq \rho_+(d_1(e,k_1))\\
&= \rho_+(d'_1(eK_i,g^{-1}hK_i))=\rho_+(d'_1(gK_i,hK_i)).
\end{align*}
Hence, taking $\rho'_-(t)=\min\{\frac{1}{3}t,\rho_-(t)\}$ and $\rho'_+(t)=\max\{3t,\rho_+(t)\}$, we have
$$\rho'_-(d'_1(gK_i,hK_i))\leq d'_2(gK_i,hK_i) \leq \rho'_+(d'_1(gK_i,hK_i)).$$ 
The functions $\rho'_-$ and $\rho'_+$ are non-decreasing, and $\lim_{t\rightarrow \infty}\rho'_\pm(t)=\infty$.
This completes the proof.
\end{proof}

\begin{Corollary}
For a finitely generated residually finite group $H$, whether or not a given box space embeds coarsely into Hilbert space is independent of the choice of finite generating set.
\end{Corollary}

\begin{Prop}\label{Subgroup}
Let $G$ be a finitely generated residually finite group with a nested sequence $\{K_i\}$ of finite index normal subgroups with trivial intersection such that the corresponding box space $\Box_{\{K_i\}}G$ embeds coarsely into Hilbert space. Then any subgroup $H$ of $G$ has a coarsely embeddable box space.
\end{Prop}
\begin{proof}
The sequence of subgroups formed by intersecting $H$ with the $K_i$ gives a sequence of finite index normal subgroups $N_i$ of $H$. Clearly the intersection of the $N_i$ is still trivial. 
Each quotient $H/N_i=H/H\cap K_i$ is isomorphic to $HK_i/K_i$, and is thus a subgroup of $G/K_i$. Thus, the box space $\Box_{\{N_i\}} H$ with the subspace metric induced from $\Box_{\{K_i\}} G$ coarsely embeds in $\Box_{\{K_i\}} G$, and hence coarsely embeds into Hilbert space. 
This metric is the metric induced on quotients of $H$ by the subspace metric on $H$ as a subgroup of $G$.
All proper left-invariant metrics on $H$ are coarsely equivalent, so by Proposition \ref{Metric}, $\Box_{\{N_i\}} H$ with a metric induced by any proper left-invariant metric on $H$ also coarsely embeds into Hilbert space.
\end{proof}

In particular, if $H$ in the above proposition is a finitely generated subgroup of $G$, then its word metric is coarsely equivalent to the metric induced on $H$ by the word metric on $G$. Hence $\Box_{\{N_i\}} H$ with the metric induced by the word metric on $H$ coarsely embeds into Hilbert space.

\begin{Corollary}\label{Free}
For any finitely generated free group $\mathbb{F}_k$, there exists a nested sequence of finite index normal subgroups with trivial intersection such that the corresponding box space embeds coarsely into Hilbert space.
\end{Corollary}
\begin{proof}
For any finite $k$, $\mathbb{F}_k$ is a finite index subgroup of $\mathbb{F}_2$, which we know has an embeddable box space by the result of \cite{AGS}.
\end{proof}

\begin{Remark}
Note that we can find an embeddable box space of a finitely generated free group $\mathbb{F}_k$ in many ways, since using the inductively defined sequence of subgroups of \cite{AGS} will also result in a coarsely embeddable box space (see Remark \ref{GenAGS}). In addition, we can view $\mathbb{F}_k$ as a subgroup of another finitely generated free group $\mathbb{F}_n$ (in many ways), and obtain the required sequence of subgroups by intersecting, as above. It would be interesting to know whether the box spaces obtained in this way are coarsely equivalent. 
\end{Remark}

\begin{Prop}\label{Overgroup}
Let $H$ be a finitely generated residually finite group with a nested sequence $\{C_i\}$ of finite index characteristic subgroups with trivial intersection such that the corresponding box space $\Box_{\{C_i\}}H$ embeds coarsely into Hilbert space. Then any group $G$ containing $H$ as a finite index normal subgroup also has a coarsely embeddable box space.
\end{Prop}
\begin{proof}
Each of the subgroups $C_i$ is normal in $G$. The box space $\Box_{\{C_i\}}G$ is coarsely equivalent to the box space $\Box_{\{C_i\}} H$, and thus embeds coarsely into Hilbert space.
\end{proof}

\begin{Example}
There exists a nested sequence of finite index normal subgroups of $SL(2,\mathbb{Z})$ with trivial intersection such that the corresponding box space embeds coarsely into Hilbert space, since $SL(2,\mathbb{Z})$ contains a finitely generated free group as a finite index normal subgroup. 
\end{Example}

\section*{Extensions}

Let $\{1\rightarrow H_i \rightarrow \Gamma_i \rightarrow G_i \rightarrow 1\}_{i\in \mathbb{N}}$ be a sequence of group extensions, where all groups involved are finite and such that the number of elements required to generate the groups $\Gamma_i$ is uniformly bounded across $i \in \mathbb{N}$. Consider the metric space $(\sqcup \Gamma_i, d_{\Gamma})$ made from the disjoint union of the $\Gamma_i$ with their word metrics, with distances between two consecutive components defined to be greater than the largest of their diameters. 
Make the disjoint union $\sqcup H_i$ into a metric space by taking the metric induced by $d_{\Gamma}$. Call this metric $d_{H}$.
Make $\sqcup G_i$ into a metric space by defining distance between two consecutive components $G_k$ and $G_{k+1}$ to be the same as the distance between the corresponding $\Gamma_k$ and $\Gamma_{k+1}$, and on each component $G_i$ taking the metric induced from $\Gamma_i$ by the quotient $\Gamma_i/H_i \cong G_i$. Call this metric $d_G$. Note that each of these spaces has bounded geometry thanks to the assumption on the generators of the $\Gamma_i$, above. 

For each $i$ let $\pi_i$ be the projection $\pi_i: \Gamma_i \longrightarrow G_i$ and choose a set-theoretic cross-section $\sigma_i: G_i \longrightarrow \Gamma_i$ such that distances to the identity are preserved. Define a map $\eta_i: \Gamma_i \times G_i \longrightarrow H_i$ by
$$\eta_i(\gamma,g)=\sigma_i(g)^{-1}\gamma \sigma_i(\pi_i(\gamma)^{-1}g).$$ We will drop the indices $i$, as this should not cause any confusion.

We will make use of the following lemma (Lemma 4.4 of \cite{DG}). Fix an index $i$.

\begin{Lemma}\label{DG}
Let $\gamma_1$ and $\gamma_2$ be elements of $\Gamma_i$, and $g$ an element of $G_i$. Then 
\begin{align*}
d_{\Gamma}(\gamma_1, \gamma_2) &\leq d_G(g,\pi(\gamma_1))+ d_G(g, \pi(\gamma_2)) + d_H(\eta(\gamma_1,g),\eta(\gamma_2,g))\\
d_H(\eta(\gamma_1,g),\eta(\gamma_2,g)) &\leq d_G(g,\pi(\gamma_1))+ d_G(g, \pi(\gamma_2)) + d_{\Gamma}(\gamma_1, \gamma_2)
\end{align*}
\end{Lemma}
We refer the reader to \cite{DG} for the proof.

\begin{Thm}\label{Ext}
Let $\{1\rightarrow H_i \rightarrow \Gamma_i \rightarrow G_i \rightarrow 1\}_{i\in \mathbb{N}}$ be as above, such that the diameters of the $\Gamma_i$ increase strictly with $i$. If the space $(\sqcup H_i,d_H)$ coarsely embeds into Hilbert space and the space $(\sqcup G_i,d_G)$ has property A, then the metric space $(\sqcup \Gamma_i, d_{\Gamma})$ embeds coarsely into Hilbert space.
\end{Thm}
\begin{proof}
We will check that the space $\sqcup \Gamma_i$ satisfies the following criterion for coarse embedding into Hilbert space given in Proposition 2.1 of \cite{DG}: if for each $\varepsilon>0$ and $R>0$ there exists a map $\varphi: \sqcup \Gamma_i \longrightarrow \ell^2(\sqcup G_i,\mathcal{H})$, $\mathcal{H}$ a real Hilbert space, such that $\|\varphi(\gamma)\|=1$ for all $\gamma \in \sqcup \Gamma_i$ and such that for all $\gamma_1, \gamma_2 \in \sqcup \Gamma_i$,
\begin{description}
\item[$(1_\Gamma)$]
if $d_{\Gamma}(\gamma_1,\gamma_2)\leq R$ then $|1-\left\langle \varphi(\gamma_1), \varphi(\gamma_2)\right\rangle| <\varepsilon$, and
\item[$(2_\Gamma)$]
$\forall \delta >0$ $\exists S>0$ such that if $d_{\Gamma}(\gamma_1,\gamma_2)\geq S$, then $|\left\langle \varphi(\gamma_1), \varphi(\gamma_2)\right\rangle| <\delta$
\end{description}
then $(\sqcup \Gamma_i,d_\Gamma)$ coarsely embeds into Hilbert space.

Let $\varepsilon>0$ and $R>0$ be given.
We know that the space $(\sqcup G_i,d_G)$ has property A, and so Proposition 2.6 of \cite{DG} tells us that there exists $\phi: \sqcup G_i \longrightarrow \ell^2(\sqcup G_i)$ and $S_G>0$ such that $\|\phi(g)\|=1$ for all $g\in \sqcup G_i$, and such that 
\begin{description}
\item[$(1_G)$]
for all $g_1, g_2 \in \sqcup G_i$, if $d_G(g_1,g_2)\leq R$ then $|1-\left\langle \phi(g_1), \phi(g_2)\right\rangle| <\frac{\varepsilon}{2}$, and
\item[$(2_G)$]
$\Supp\phi(g) \subset B_{S_G}(g)$ for all $g\in \sqcup G_i$.
\end{description}
We will view $\phi$ as a function on $\sqcup G_i \times \sqcup G_i$.
Since $\sqcup H_i$ is coarsely embeddable, there exists according to Proposition 2.1 of \cite{DG} a Hilbert space-valued map $\psi: \sqcup H_i \longrightarrow \mathcal{H}$ such that $\|\psi(h)\|=1$ for all $h \in \sqcup H_i$ and such that for all $h_1, h_2 \in \sqcup H_i$,
\begin{description}
\item[$(1_H)$]
if $d_{H}(h_1,h_2)\leq 2S_G+R$ then $|1-\left\langle \psi(h_1), \psi(h_2)\right\rangle| <\frac{\varepsilon}{2}$, and
\item[$(2_H)$]
$\forall \delta >0$ $\exists S_H>0$ such that if $d_{H}(h_1,h_2)\geq S_H$, then $|\left\langle \psi(h_1), \psi(h_2)\right\rangle| <\delta$.
\end{description}

We will now define the map $\varphi: \sqcup \Gamma_i \longrightarrow \ell^2(\sqcup G_i, \mathcal{H})$. Let $N_R\in \mathbb{N}$ be such that for $i\geq N_R$, the distance between $\Gamma_i$ and $\Gamma_{i+1}$ is greater than $R$. For $\gamma \in \Gamma_i$ with $i\leq N_R$, let $\varphi(\gamma)$ be given by $\varphi(\gamma)(\pi(e_1))=\psi(e_1)$ and $\varphi(\gamma)(g)=0$ for $g\neq \pi(e_1)$, where $e_1$ denotes the identity element of $\Gamma_1$. For $\gamma\in \Gamma_i$ with $i> N_R$, define
$$\varphi (\gamma)(g)=\phi(\pi(\gamma), g)\psi(\eta(\gamma,g))$$
for $g\in G_i$ and $$\varphi(\gamma)(g)=\phi(\pi(\gamma), g)\psi(e_i)$$ for $g\notin G_i$, where $e_i$ is the identity element of $H_i \triangleleft G_i$.

We will first check that $\|\varphi(\gamma)\|=1$ for all $\gamma \in \sqcup \Gamma_i$. For $\gamma\in \Gamma_i$ with $i\leq N_R$, we have 
\begin{equation*}
\|\varphi(\gamma)\|=\|\psi(e_1)\|=1
\end{equation*}
and for $\gamma\in \Gamma_i$ with $i> N_R$, we have
\begin{align*}
\|\varphi(\gamma)\|^2&= \sum_{g\in G_i} |\phi(\pi(\gamma),g)|^2\cdot \|\psi(\eta(\gamma,g))\|^2 
+ \sum_{g\notin G_i} |\phi(\pi(\gamma),g)|^2 \cdot \|\psi(e_i)\|^2\\
&= \sum_{g\in G_i} |\phi(\pi(\gamma),g)|^2 + \sum_{g\notin G_i} |\phi(\pi(\gamma),g)|^2\\
&= \|\phi(\pi(\gamma))\|^2= 1.
\end{align*}

Let us now check the remaining two conditions $(1_{\Gamma})$ and $(2_{\Gamma})$.

Take $\gamma_1$ and $\gamma_2$ in $\sqcup \Gamma_i$ such that $d_{\Gamma}(\gamma_1,\gamma_2)\leq R$. Then if both $\gamma_1$ and $\gamma_2$ lie in $\sqcup_{i=1}^{N_R}\Gamma_i$, we have $|1-\left\langle \varphi(\gamma_1), \varphi(\gamma_2)\right\rangle|=|1-\|\psi(e_1)\|^2|=|1-1|=0$. 
The only other possibility is that they both lie in the same component, say $\Gamma_i$, for $i> N_R$. In this case we have
\begin{align*}
|1-\left\langle \varphi(\gamma_1), \varphi(\gamma_2)\right\rangle|= &\bigg|1-\sum_{g\in G_i}\phi(\pi(\gamma_1),g)\phi(\pi(\gamma_2),g)\left\langle \psi(\eta(\gamma_1,g)),\psi(\eta(\gamma_2,g))\right\rangle\\ 
&- \sum_{g\notin G_i}\phi(\pi(\gamma_1),g)\phi(\pi(\gamma_2),g)\left\langle \psi(e_i),\psi(e_i) \right\rangle \bigg| \\
= &\bigg| 1 -\sum_{g\in G_i}\phi(\pi(\gamma_1),g)\phi(\pi(\gamma_2),g)\left\langle \psi(\eta(\gamma_1,g)),\psi(\eta(\gamma_2,g))\right\rangle\\ 
&+ \sum_{g \in G_i}\phi(\pi(\gamma_1),g)\phi(\pi(\gamma_2),g)- \sum_{g \in G_i}\phi(\pi(\gamma_1),g)\phi(\pi(\gamma_2),g)\\
&- \sum_{g\notin G_i}\phi(\pi(\gamma_1),g)\phi(\pi(\gamma_2),g) \bigg| \\
\leq &\bigg| \sum_{g\in G_i} \phi(\pi(\gamma_1),g)\phi(\pi(\gamma_2),g)(1-\left\langle \psi(\eta(\gamma_1,g)),\psi(\eta(\gamma_2,g))\right\rangle)\bigg| \\
+ &\bigg|1- \left\langle\phi(\pi(\gamma_1)),\phi(\pi(\gamma_2))\right\rangle\bigg|
\end{align*}
The quotient map $\pi$ is contractive and hence we have that $d_G(\pi(\gamma_1),\pi(\gamma_2))\leq R$, so $(1_G)$ tells us that the second term is bounded by $\frac{\varepsilon}{2}$. By $(2_G)$, the sum in the first term ranges over $g \in B_{S_G}(\pi(\gamma_1))\cap B_{S_G}(\pi(\gamma_2))$, and thus the first term can be bounded by
$$\sup\{|1-\left\langle \psi(\eta(\gamma_1,g)),\psi(\eta(\gamma_2,g))\right\rangle|: g\in B_{S_G}(\pi(\gamma_1))\cap B_{S_G}(\pi(\gamma_2))\}.$$
By Lemma \ref{DG} above, for $g\in B_{S_G}(\pi(\gamma_1))\cap B_{S_G}(\pi(\gamma_2))$ we have $$d_H(\eta(\gamma_1,g),\eta(\gamma_2,g))\leq 2S_G + R$$ and so by $(1_H)$, the supremum above is bounded by $\frac{\varepsilon}{2}$. This completes the proof of $(1_{\Gamma})$.

For $(2_{\Gamma})$, fix $\delta >0$. Take the $S_H$ corresponding to $\frac{\delta}{3}$ as in $(2_H)$. The required $S$ will be given by $3S_G + 3S_H+M_{\Gamma}$, where $$M_{\Gamma}= \max\{d_{\Gamma}(\gamma_1,\gamma_2): \gamma_1,\gamma_2 \in \sqcup_{i=1}^{N_R} \Gamma_i \}.$$ 

If $d_{\Gamma}(\gamma_1,\gamma_2)\geq 3S_G + 3S_H+M_\Gamma$, then at most one of $\gamma_1,\gamma_2$ lies in $\sqcup_{i=1}^{N_R} \Gamma_i$. If one of them does, without loss of generality $\gamma_1 \in \sqcup_{i=1}^{N_R} \Gamma_i$, and $\gamma_2\in \Gamma_j$ for some $j>N_R$. We then have
\begin{align*}
|\left\langle \varphi(\gamma_1), \varphi(\gamma_2) \right\rangle| &=|\phi(\pi(\gamma_2),\pi(e_1))\left\langle \psi(e_1),\psi(e_j) \right\rangle|.
\end{align*}
Now notice that $d_H(e_1,e_j)=d_{\Gamma}(e_1,e_j)\geq d_{\Gamma}(\gamma_1,e_j) \geq d_{\Gamma}(\gamma_1,\gamma_2)- d_{\Gamma}(\gamma_2,e_j)$, and
\begin{align*}
2d_{\Gamma}(\gamma_2,e_j)&\leq \diam(\Gamma_j) +d_{\Gamma}(\gamma_2,e_j)\\
&\leq d_{\Gamma}(e_{j-1},e_j) +d_{\Gamma}(\gamma_2,e_j) \\
&\leq d_{\Gamma}(\gamma_1,\gamma_2),
\end{align*}
so $d_H(e_1,e_j)$ is greater than or equal to $d_{\Gamma}(\gamma_1,\gamma_2) - \frac{1}{2}d_{\Gamma}(\gamma_1,\gamma_2)\geq S_H$. We thus have $|\left\langle \varphi(\gamma_1), \varphi(\gamma_2) \right\rangle| < \delta$ by $(2_H)$.

If neither $\gamma_1$ nor $\gamma_2$ lie in $\sqcup_{i=1}^{N_R} \Gamma_i$, then there are two possibilities. If $\gamma_1 \in \Gamma_i$ and $\gamma_2\in \Gamma_j$ for $N_R<i< j$, then we have
\begin{align*}
|\left\langle \varphi(\gamma_1), \varphi(\gamma_2) \right\rangle| &\leq \bigg| \sum_{g\in G_i} \phi(\pi(\gamma_1),g)\phi(\pi(\gamma_2),g)\left\langle \psi(\eta(\gamma_1,g)),\psi(e_j)\right\rangle\bigg|\\ &+\bigg|\sum_{g\in G_j}\phi(\pi(\gamma_1),g)\phi(\pi(\gamma_2),g)\left\langle \psi(e_i),\psi(\eta(\gamma_2,g))\right\rangle\bigg|\\ &+ \bigg|\sum_{g\notin G_i, G_j}\phi(\pi(\gamma_1),g)\phi(\pi(\gamma_2),g)\left\langle \psi(e_i),\psi(e_j)\right\rangle\bigg|.
\end{align*}
The distances between each of the pairs of points $(\eta(\gamma_1,g),e_j)$, $(e_i,\eta(\gamma_2,g))$ and $(e_i,e_j)$ are greater than $d_{\Gamma}(e_i,e_j)$ and thus greater than $\frac{1}{3}d_{\Gamma}(\gamma_1,\gamma_2) \geq S_H$ 
(since $d_\Gamma(\gamma_1,\gamma_2)\leq d_{\Gamma}(\gamma_1,e_i)+d_{\Gamma}(e_i,e_j)+d_{\Gamma}(e_j,\gamma_2)\leq 3d_{\Gamma}(e_i,e_j)$),
so by $(2_H)$ we have $|\left\langle \varphi(\gamma_1), \varphi(\gamma_2) \right\rangle| < \delta$ as required.

The second possibility is that $\gamma_1$ and $\gamma_2$ lie in the same component $\Gamma_i$, for $i>N_R$. Then we have
\begin{align*}
|\left\langle \varphi(\gamma_1), \varphi(\gamma_2) \right\rangle| &\leq \bigg| \sum_{g\in G_i} \phi(\pi(\gamma_1),g)\phi(\pi(\gamma_2),g)\left\langle \psi(\eta(\gamma_1,g)),\psi(\eta(\gamma_2,g))\right\rangle\bigg|\\
&+ \bigg|\sum_{g\notin G_i}\phi(\pi(\gamma_1),g)\phi(\pi(\gamma_2),g)\left\langle \psi(e_i),\psi(e_i)\right\rangle\bigg|
\end{align*}
The first term is bounded above by 
$$\sup\{|\left\langle \psi(\eta(\gamma_1,g)),\psi(\eta(\gamma_2,g))\right\rangle|: g\in B_{S_G}(\pi(\gamma_1))\cap B_{S_G}(\pi(\gamma_2))\}.$$
So, from Lemma \ref{DG} above, we deduce that $d_H(\eta(\gamma_1,g),\eta(\gamma_2,g))\geq d_{\Gamma}(\gamma_1,\gamma_2)-d_G(g,\pi(\gamma_1))-d_G(g,\pi(\gamma_2))\geq S_H$ for any $g\in B_{S_G}(\pi(\gamma_1))\cap B_{S_G}(\pi(\gamma_2))$ and thus the first term is bounded by $\frac{\delta}{3}$. 

Now $B_{S_G}(\pi(\gamma_1))\cap B_{S_G}(\pi(\gamma_2)) \subset G_i$, since for any $\gamma$ in $\Gamma_i$ and $g$ in $G_k$, $k\neq i$, we have
$$d_G(g,\pi(\gamma))\geq d_G(e_k,e_i)\geq \diam(\Gamma_i)\geq d_{\Gamma}(\gamma_1,\gamma_2) \geq S_G.$$
Thus, the set over which the second sum is taken is empty. This completes the proof of $(2_{\Gamma})$, and hence the theorem.

\end{proof}

Note that we do not require the diameters of the $H_i$ or the $G_i$ to increase. 

\section*{Applications}
In this section, we apply Theorem \ref{Ext} in two different ways. First, we start with a residually finite group $\Gamma$, which is an extension of $H$ by $G$, and some sequence of nested finite index normal subgroups of $\Gamma$ with trivial intersection. We give a sufficient condition for the corresponding box space to embed coarsely into Hilbert space, in terms of the groups $G$ and $H$. 

However, in practice it may be difficult to check whether these conditions hold for a given group. We therefore give a concrete application of Theorem \ref{Ext} for semidirect products, where we can build a sequence of nested finite index normal subgroups of the semidirect product out of such sequences for the factors.

Consider first a sequence of extensions arising from an extension $1\longrightarrow H \longrightarrow \Gamma \longrightarrow G \longrightarrow 1$ in the following way. Suppose $\Gamma$ is residually finite, and let $\{K_i\}$ be a sequence of nested finite index normal subgroups of $\Gamma$ with trivial intersection. Then $H\cap K_i$ is such a sequence for $H$, and we have the following sequence of extensions:
$$1 \longrightarrow H/H\cap K_i \longrightarrow \Gamma/K_i \longrightarrow \Gamma/HK_i \longrightarrow 1,$$
where the groups $\Gamma/HK_i$ can be seen as finite quotients of $G$ by $L_i:=HK_i/H$. We can then apply our theorem to conclude that the box space $\Box_{\{K_i\}} \Gamma$ coarsely embeds into Hilbert space if the space $\sqcup G/L_i$ has property A and the box space $\Box_{\{H\cap K_i\}} H$ with the induced metric coarsely embeds into Hilbert space. Note that in general, $G$ will not be residually finite.

Define $L$ to be the intersection of the $L_i$.

\begin{Prop}
Let $1\longrightarrow H \longrightarrow \Gamma \longrightarrow G \longrightarrow 1$ be an extension as above. If $G/L$ is amenable and the box space $\Box_{\{H\cap K_i\}} H$ with the induced metric coarsely embeds into Hilbert space, then $\Box_{\{K_i\}} \Gamma$ embeds coarsely into Hilbert space.
\end{Prop}

\begin{proof}
For each $i$, we have $G/L_i \cong (G/L)/(L_i/L)$. Note that the intersection $\cap L_i/L$ is trivial, and that the quotient $G/L$ is amenable. Hence, by Guentner's result, the box space $\Box_{\{L_i/L\}} G/L$ has property A and hence the space $\sqcup G/L_i$ does too. We can now apply Theorem \ref{Ext} to conclude that $\Box_{\{K_i\}} \Gamma$ embeds coarsely into Hilbert space.
\end{proof}

We can now state a sufficient condition for a residually finite group to have the Haagerup property.

\begin{Corollary}
Let $\Gamma$ be a finitely generated residually finite group which is an extension of a group $H$ by a group $G$. Suppose there exists a sequence $\{K_i\}$ of nested normal finite index subgroups of $\Gamma$ such that the intersection $\cap K_i$ is trivial and $\Box_{\{H\cap K_i\}} H$ with the induced subspace metric coarsely embeds into Hilbert space. Let $L$ denote the intersection of the subgroups $HK_i/H$ of $G$. Then $\Gamma$ has the Haagerup property if $G/L$ is amenable.
\end{Corollary}

Let us remark here that the Haagerup property is known to be preserved under extensions with amenable quotients (see \cite{CJV}).

Since it may be difficult to check that the space $\Box_{\{H\cap K_i\}} H$ is embeddable, we now give a more concrete application for semidirect products.

\begin{Thm}\label{Semi}
Let $\Gamma$ be the semidirect product $H\rtimes G$ of two finitely generated residually finite groups $H$ and $G$ such that there is a nested sequence of finite index characteristic subgroups $\{N_i\}$ of $H$ with $\cap N_i = 1$, such that $\Box_{\{N_i\}}H$ embeds coarsely into Hilbert space, and G is amenable. Then $\Gamma$ has an embeddable box space.
\end{Thm}

\begin{proof}
Enumerate the non-trivial elements of $\Gamma$, so that $\Gamma = \{e, \gamma_1, \gamma_2, \gamma_3, ... \}$. For each $\gamma_i$, we will find a normal finite index subgroup $K_i$ of $\Gamma$ such that the image of $\gamma_i$ is non-trivial in $\Gamma/K_i$ and such that the $K_i$ are nested. We will do this inductively. We essentially prove that $\Gamma$ is residually finite while making sure that the subgroups $K_i$ are built from the subgroups of $H$ and $G$ in a particular way.

First, take $\gamma_1= (x,a)$, where $x\in H$ and $a\in G$. If $a$ is non-trivial, then there is some finite quotient $Q$ of $G$ in which the image of $a$ is still non-trivial. In this case, take $K_1$ to be the kernel of the surjection $\Gamma \longrightarrow Q$. If $a$ is trivial, then $x$ is non-trivial, and so we can take one of the characteristic subgroups $N_j$ of H such that the image of $x$ is non-trivial in the quotient $H/N_j$. Since $N_j$ is characteristic in $H$, it is normal in $\Gamma$ and so we have a quotient of $\Gamma$ which is isomorphic to $H/N_j \rtimes G$, in which the image of $\gamma_1$ is non-trivial. To get a finite quotient, we take the subgroup $A$ of $G$ which acts trivially on $H/N_j$. $A$ is normal in $H/N_j \rtimes G$, and the quotient $H/N_j \rtimes G/A$ is finite because each element of $G/A$ now acts as a non-trivial automorphism of $H/N_j$, and so $G/A$ is a subgroup of $\Aut(H/N_j)$, which is finite since $H/N_j$ is. Let $K_1$ be the kernel of the homomorphism $\Gamma \longrightarrow H/N_j \rtimes G/A$.

Now, suppose we have defined $K_1, K_2, ...$ up to $K_{i-1}$. Given $\gamma_i$, we want to find a normal finite index subgroup $K_i$ of $\Gamma$ such that the image of $\gamma_i$ is non-trivial in $\Gamma/K_i$ and such that $K_i \leq K_{i-1}$. Let $\gamma_i$ be given by $(h,b)$, where $h\in H$ and $b\in G$. Take the characteristic subgroup $N_k$ of $H$ such that $K_{i-1}\cap H \geq N_k$ and, additionally, such that the image of $h$ is non-trivial in $H/N_k$ if $h$ is non-trivial. Consider the quotient $H/N_k \rtimes G$. Now take the subgroup $B$ of $G$ which acts trivially on $H/N_k$, and take the intersection with $K_{i-1}\cap G$. Call this subgroup $K$. Now $K$ still acts trivially on $H/N_k$, and is thus a normal subgroup of $H/N_k \rtimes G$. The quotient $H/N_k \rtimes G/K$ is finite (since $K_{i-1}\cap G$ and $B$ are both of finite index in $G$, and hence so is their intersection), and $\gamma_i$ has a non-trivial image in this quotient. Thus, we can define $K_i$ to be the kernel of the map $\Gamma \longrightarrow H/N_k \rtimes G/K$, which lies inside $K_{i-1}$ by construction. 

We can now apply Theorem \ref{Ext}. Note that we have obtained a sequence of nested normal finite index subgroups $\{K_i\}$ such that each quotient $\Gamma/K_i$ is an extension of a quotient $H/N_n$ by a finite quotient of $G$. The intersection of all the $K_i$ is trivial. The disjoint union of the quotients of $H$ coarsely embeds into Hilbert space, and the disjoint union of the quotients of $G$ has property A, since $G$ is amenable. The conditions of Theorem \ref{Ext} are thus satisfied, so the box space $\Box_{\{K_i\}} \Gamma$ coarsely embeds into Hilbert space. 
\end{proof}

\begin{Remark}
We can apply the above to show that semidirect products of finitely generated free groups by residually finite amenable groups have an embeddable box space. This provides a new class of examples of spaces with bounded geometry which embed coarsely into Hilbert space but do not have property A, generalising the example of Arzhantseva, Guentner and Spakula.
\end{Remark}

\begin{Corollary}
Let $\Gamma$ be an extension of a finitely generated free group by a (finite or infinite) cyclic group. Then $\Gamma$ has an embeddable box space.
\end{Corollary}

\begin{proof}
Suppose first that $\Gamma$ is an extension of a finitely generated free group $\mathbb{F}_n$ by a finite cyclic group $C$. Take the sequence of subgroups of $\mathbb{F}_n$ defined inductively as in \cite{AGS}, $\{N_i\}$. These subgroups are characteristic in $\mathbb{F}_n$, and hence we can apply Proposition \ref{Overgroup}, which gives the required result.

Suppose now that $\Gamma$ is an extension of a finitely generated free group $\mathbb{F}_n$ by $\mathbb{Z}$. This extension splits, and so $\Gamma$ is isomorphic to $\mathbb{F}_n \rtimes \mathbb{Z}$. 
We can now apply Theorem \ref{Semi} to conclude the proof.
\end{proof}

\end{document}